\documentclass[MSNbibl,number,citesort,rotating,dvips]{arxbj}
\usepackage{dcolumn}
\usepackage{graphicx}

\aid{0}
\volume{19}
\issue{1}
\pubyear{2013}
\firstpage{344}
\lastpage{361}
\doi{10.3150/11-BEJ388}

\makeatletter
\def\bsuffix #1{#1}
\newremark{remark}{Remark}
\newremark{example}{Example}
\newtheorem{theorem}{Theorem}
\newtheorem{lemma}{Lemma}
\newcolumntype{d}[1]{D{.}{.}{#1}}
\makeatother

\begin{document}
\begin{frontmatter}

\title{A quantile regression estimator for~censored~data}
\runtitle{Censored quantile regression}

\begin{aug}
%
\author[1]{\fnms{Chenlei} \snm{Leng}\corref{}\thanksref{1}\ead[label=e1]{stalc@nus.edu.sg}}%
\and
\author[2]{\fnms{Xingwei} \snm{Tong}\thanksref{2}\ead[label=e2]{xweitong@bnu.edu.cn}}
\runauthor{C. Leng and X. Tong} 
\address[1]{Department of Statistics and Applied Probability,
National University of Singapore, 6 Science Drive 2, SG 117546,
Republic of Singapore. \printead{e1}}
\address[2]{School of Mathematical Sciences, Beijing Normal University,
Beijing 100875,
China.\\
\printead{e2}}
\end{aug}

\received{\smonth{3} \syear{2010}}
\revised{\smonth{5} \syear{2011}}

%
\begin{abstract}
We propose a censored quantile regression estimator motivated by
unbiased estimating equations. Under the usual conditional independence
assumption of the survival time and the censoring time given the
covariates, we show that the proposed estimator is consistent and
asymptotically normal. We develop an efficient computational
algorithm which uses existing quantile regression code. As a result,
bootstrap-type inference can be efficiently implemented.
We illustrate the finite-sample
performance of the proposed method
by simulation studies and analysis of a survival data set.
\end{abstract}

%
\begin{keyword}
\kwd{accelerated failure time model}
\kwd{censored quantile regression}
\kwd{Kaplan--Meier estimate}
\kwd{quantile regression}
\end{keyword}

\end{frontmatter}

\section{Introduction}\label{sec1}
Censored data arise frequently in biomedical, psychological, social
studies and many other applied fields (Kalbfleisch and Prentice~\cite{Kal2002}).
Analysis of such data is complicated by censoring, where an object's
time-to-death or other end-point of interest is known to occur only in
a certain period of time.

Denote $T$ as the survival time and $C\leq\mathcal{ T}_0$ as the
censoring time, where $\mathcal{ T}_0$ is the largest follow-up study
time. The typical censored data set consists of independent
observations $(Y_i,\delta_i,Z_i)$, $i=1,\ldots, n$, where
$Y_i=\operatorname{min}(T_i,C_i)$ is the observed failure time;
$\delta_i=I(T_i \le C_i)$ is the censoring indicator; and $Z_i$
is a $p$-dimensional covariate vector including an intercept.
The accelerated failure time model (AFT), specified as $T_i=\beta^T
Z_i+\varepsilon_i$ with $\varepsilon_i, i=1,\ldots ,n$, following a
common distribution independently, was studied in a number of papers
(Jin \textit{et al.}~\cite{Jinetal2003}; Zeng and Lin~\cite{ZengLin2007}). However, the assumption made
for the AFT model precludes error heteroscedasiticity and will yield
biased results when it is inappropriate.

Quantile regression introduced by Koenker and Bassett~\cite{Koenker1978} has
become an increasingly important tool in statistical analysis.
Contrary to the usual model for the conditional mean, it provides
distributional information on the dependence of $T$ on $Z$.
A~comprehensive review can be found in Koenker~\cite{Koenker2005}. The usefulness
of quantile regression in survival analysis was discussed by Koenker
and Geiling~\cite{Koenker2001}.
The $\tau$th conditional quantile function of the dependent variable
$T$ given covariates $Z$, $Q_T(\tau|Z)$, is defined as
$Q_T(\tau|Z)=\operatorname{inf}\{v\dvt F_0(v|Z) \ge\tau\}$, where $F_0$ is the
cumulative conditional distribution function of $T$ given $Z$.
Correspondingly, a quantile regression model for $Q_T(\tau|Z)$ with
$\tau\in(0,1)$ can be denoted as
%
\begin{equation}
Q_T(\tau|Z)=\beta_\tau'Z. \label{qr}
\end{equation}
Note that the AFT
model is a special case of this model when $\beta_1$, the coefficient
corresponding to the intercept, is quantile dependent, but the other
coefficients in $\beta$ are quantile independent. Another
example is the location-scale model $T_i=\alpha_1^T Z_i+(\alpha_2^T
Z_i)\varepsilon_i$ in which $\beta=\alpha_1+\alpha_2
F^{-1}_\varepsilon(\tau)$ when $\varepsilon_i$ is independent of
$Z_i$, and $F_\varepsilon$ is the distribution function of
$\varepsilon_i$.

When data are subject to censoring, statistical estimation and
inference for quantile regression is more involved. Indeed, a naive
procedure which completely ignores censoring may give highly biased
estimates (Koenker~\cite{Koenker2005}). The situation is more complicated if censoring
time depends on the covariates.
In Section~\ref{sec5}, we present the Colorado
Plateau uranium miners cohort data (Lubin \textit{et al.}~\cite{L1994}, Langholz and
Goldstein~\cite{Lang1996}), where the major interest of this study is to assess the
effect of smoking and radon exposure on the rate of median death time of
lung cancer. It is found that the censoring time is highly correlated
with the
covariates. Ignoring this dependence may yield biased estimates; see
the Numerical Study section for examples.

Powell~\cite{Powell1984,Powell1986} first
studied censored quantile regression with fixed censoring. For random
censoring, Ying, Jung and Wei~\cite{Yingetal1995} (YJW) proposed a semiparametric
median regression model. Despite the simplicity of the method in YJW,
this procedure requires the unconditional independence of the survival
time and censoring time. This assumption is often restrictive as
conditional independence, given the covariates, is more natural
(Kalbfleisch and Prentice~\cite{Kal2002}).
In addition, the estimating equation approach proposed in YJW involves
solving non-monotone discrete equations, creating difficulty for optimization.
As a consequence, inferential procedures such as the resampling
approach in Jin, Ying and Wei~\cite{Jinetal2001}, or the bootstrap method, can be
prohibitive computationally. See also Leon, Cai and Wei~\cite{Leon2009} for a
generalization of this method to partly linear models.

Relaxing the independence condition to conditional independence,
Portnoy~\cite{Portnoy2003} and Neocleous, Vanden Brandan
and Portnoy~\cite{Neo2006} developed a novel estimating approach motivated by
the classical
Kaplan--Meier estimator in the one sample analysis. Using
the martingale representation, Peng and Huang~\cite{PengHuang2008} studied another
censored quantile regression estimator motivated by the Nelson--Aalen
estimator. However, a major
shortcoming of Portnoy and Peng and Huang's methods is that a global
linear assumption has to be made, even for estimating the quantile
coefficient at a single quantile. To relax this
assumption, Wang and Wang~\cite{WangWang2009} recently proposed an innovative
redistribution of mass idea, which employs local weighting.

Motivated by the unbiased estimating equation for the quantile
regression (Ying \textit{et al.}~\cite{Yingetal1995}), we propose a new quantile
regression estimator. Under the usual conditional independence
assumption of $T$ and $C$ given $Z$, we show that the proposed
estimator is consistent and asymptotically normal. We develop an
efficient algorithm which utilizes existing quantile regression code
for estimation. The efficient code enables us to use the bootstrap
procedure for statistical inference. Our method provides an alternative
to Wang and Wang's locally weighted censored quantile regression.
However, the framework proposed by Ying \textit{et al.}, and used by us,
may be conceptually simpler.

The rest of the paper is organized as follows. Section~\ref{sec2} discusses the
new estimator and the fast computing algorithm. Section~\ref{sec3} provides the
theoretical results of the new estimator. Some numerical studies are
presented in Section~\ref{sec4}. A data analysis is provided in Section~\ref{sec5}. Some
concluding remarks are given in Section~\ref{sec6}. All the proofs are relegated
to the \hyperref[app]{Appendix}.
When no confusion arises, the dependence of $\beta$ on $\tau$ is suppressed.

\section{Censored quantile regression}\label{sec2}
To estimate $\beta_\tau$ in (\ref{qr}) for the $\tau$th quantile,
we propose to solve the following estimating equation:
%
\begin{equation}
M_n(\beta)=\sum_{i=1}^n Z_i\biggl[\frac{I(Y_i-\beta'Z_i\geq0)}{\hat
{G}(\beta' Z_i|Z_i)} -(1-\tau) \biggr] \approx0, \label{eqyjw}
\end{equation}
where $I(\cdot)$ is the indicator function, and $\hat{G}$ is the
Kaplan--Meier estimate for $G_0(\cdot|Z_i)$, the conditional
survival function of the censoring variable $C$ given the covariates.
The estimating equation in (\ref{eqyjw}) is motivated by the fact
that $E[I(Y_i-\beta'Z_i \ge0)|Z_i]=(1-\tau)G_0(\beta'Z_i|Z_i)$
using the conditional independence of $T_i$ and $C_i$ given $Z_i$.
YJW assumes that $G_0(\beta'Z_i|Z_i)=G_0(\beta'Z_i)$, and our formulation
is an extension of YJW's median regression to
quantile regression by allowing $G_0$ to
depend on $Z$. To solve (\ref{eqyjw}), we need to find its root. Ying
\textit{et al.}~\cite{Yingetal1995}
proposed to minimize $\|M_n(\beta) \|$, a discrete and non-monotone
function. Computational complication naturally arises.
Ying \textit{et al.} proposed to use the simulated
annealing algorithm, or simply the bisection algorithm, to solve the
estimating equation, which is computationally demanding.
Another complication arises in statistical inference.
Since the sampling distribution of the solution
involves the unknown density functions of the data, resampling-based
approaches are effective tools for conducting statistical inference (Jin,
Lin, Wei and Ying~\cite{Jinetal2003}). However, inference procedures via these
methods would be computationally even more intensive than point
estimation, due to the lack of an efficient algorithm.

We start by proposing a new algorithm to solve (\ref{eqyjw}). First
note that we can write the estimating equation in (\ref{eqyjw}) as
\[
M_n(\beta)=\sum_{i=1}^n \frac{Z_i}{G_i}[I(Y_i-\beta'Z_i \ge
0)-(1-\tau) ]-(1-\tau) \sum_{i=1}^n \frac{Z_i}{G_i}(G_i-1),
\]
where we write $\hat{G}(\beta' Z_i|Z_i)$ for a preliminary estimate of
$\beta$ as $G_i$ for brevity.
In practice, we set $I(Y_i-\beta^TZ_i \ge0)/G_i=0$ if $G_i=0$ as
recommended by Ying \textit{et al.}~\cite{Yingetal1995}.
The solution to this function is the minimizer of the following
linear programming problem:
%
\begin{equation}S_n(\beta)= \sum_{i=1}^n G_i^{-1} \bigl\{
\rho_\tau(Y_i-\beta^T Z_i)
+\rho_\tau\bigl(Y^*_i- \beta^T Z_i(G_i-1) \bigr) \bigr\},
\label{eqqr}
\end{equation}
where $\rho_\tau(s)=s[I(s \ge0)-(1-\tau)]$
is the check loss function used in quantile regression,
and $Y^*_i$ is a small constant less than
$-|\beta^TZ_i(G_i-1)|$ for all $\beta$'s in a compact space.
In this paper, we set
$Y^*_i=Y^*=\min\{Y_i \}-A$ with $A=200$. This new formulation suggests that
the fast quantile regression code (Portnoy and Koenker~\cite{Portnoy1997}, Koenker
\cite{Koenker2005}) can be directly used to
solve the censored quantile regression defined by (\ref{eqyjw}). We
note that
Yin and Cai~\cite{Yin2005} used the Nelder--Mead simplex algorithm when
$T$ and $C$ are independent. The simplex
algorithm is generally much slower than the interior point algorithm specially
developed for quantile regression (Portnoy and Koenker~\cite{Portnoy1997}).

For estimating the weighting function $G_0(\beta'Z_i|Z_i)$,
we propose to use the
local Kaplan--Meier estimator $\hat{G}(\beta'Z_i|Z_i)$. To be specific,
$G_0(\cdot|Z_i)$ is estimated by
%
\begin{equation}\hat{G}(t|z)=\prod_{j=1}^n\biggl[1-\frac
{B_{nj}(z)}{\sum_{k=1}^n
I(Y_k\geq Y_j) B_{nk}(z)}\biggr]^{I(Y_j\leq t,\delta_j=0)},
\label{eqEstConC}
\end{equation}
where $B_{nj}(z)$ is a sequence of
non-negative weights adding up to 1. When $B_{nj}(z)=1/n$ for all
$j$, $\hat{G_0}(t|z)$ reduces to the classical Kaplan--Meier
estimator of the survival function in the one-sample case.
Following the idea
of Wang and Wang~\cite{WangWang2009}, we use
%
\begin{equation}
B_{nj}(z)=K\biggl(\frac{z-z_j}{h_n}\biggr)\Biggl[\sum_{k=1}^n
K\biggl(\frac{z-z_k}{h_n}\biggr)\Biggr]^{-1}, \label{eqEstB}
\end{equation}
where
$K(\cdot)$ is a density function, and $h_n>0$ is the
bandwidth. This is the familiar kernel estimator for the survival function
discussed, for example, in Gonzalez-Manteiga and Cadarso-Suarez~\cite{G1994}.
When $Z$ is multi-dimensional, we can use the product kernel. For
example, in the bivariate case, we can use
$K(z_1,z_2)=K_1(z_1)K_2(z_2)$ where $K_1(\cdot)$
and $K_2(\cdot)$ are both univariate kernel functions. In this
article, we
use the bi-quadratic kernel, defined as
$K(s)=\frac{15}{16}(1-s^2)^2I(|s|\le1)$, which is also used
by Wang and Wang~\cite{WangWang2009}.
Alternatively, we may use a multivariate density function, for example,
from the multivariate normal distribution (Fan and Gijbels~\cite{FanGijbels1996}).

Since $\hat{G}(\beta'Z_i|Z_i)$ depends on the unknown parameter
$\beta
$, we propose the following iterative algorithm between solving
for $\hat\beta$ and $\hat{G}(\beta'Z_i|Z_i)$, while the other one
is fixed:
\begin{enumerate}[1.]
\item[1.] Given an initial estimate of $\beta$ denoted as $\beta
^{(0)}$, set $k=0$.
\item[2.] Estimate $G_0(Z_i'\beta^{(k)}|Z_i)$ as $G_i$
using the local Kaplan--Meier estimator. Minimize $S_n(\beta)$ in
(\ref
{eqqr}) to obtain $\beta^{(k+1)}$.
\item[3.] Set $k \leftarrow k+1$. Go to Step 2 until a convergence
criterion is met.
\end{enumerate}
For the initial estimate, we use a similar method as in
(Yin and Cai~\cite{Yin2005}) to solve the following monotone estimating function:
\[
\sum_{i=1}^n \frac{\delta_i}{\hat{G}(Y_i|Z_i)}Z_i[I(Y_i-\beta
^TZ_i\ge
0)-(1-\tau)],
\]
where $\hat{G}(Y_i|Z_i)$ is the local Kaplan--Meier estimator.
This estimator can be seen as the inverse probability weighted quantile
regression function (Bang and Tsiatis~\cite{bang2002}). Similarly to (Yin and Cai
\cite{Yin2005}), consistency of $\beta^{(0)}$ can be shown, and convergence of
the solution series $\{\beta^{(k)} \}$ follows by the method of
induction. Note that, although the initial estimate is also reasonable,
its efficiency is adversely affected by the fact that only non-censored
observations are used.

\begin{remark}\label{rem1} We note that our estimator requires estimating
$G_0(\beta'Z_i|Z_i)$. 
In comparison, Wang and Wang~\cite{WangWang2009} estimated the conditional
cumulative distribution function of the survival time $T_i$, evaluated
at $C_i$, given the covariate $Z_i$. Both estimators use local
Kaplan--Meier estimation. Since both methods are estimation equation
based approaches, which one is more efficient is likely
dependent on the particular problem under analysis. We observe
empirically that the proposed method performs satisfactorily even if
the censoring rate is reasonably low. Of course, if the censoring rate
is very low, the proposed method ultimately suffers due to the low
sample size used for the local Kaplan--Meier estimate.

We briefly discuss the computation issue before presenting the
asymptotic results. The estimation problem in (\ref{eqqr}) is
essentially weighted quantile regression after the weights
$G_0(Z_i'\beta|Z_i)$ are estimated using the local Kaplan--Meier
method. This method can be easily implemented by extending the
Kaplan--Meier estimate for the survival function of $C$. In particular,
we augment observations $(Y_i^*, Z_i(G_i-1))$ with weights $1/G_i$
to the existing data set $(Y_i, Z_i)$ with weight $1/G_i$. We then
apply function \textit{rq} in R library \textit{quantreg} on the augmented
data set using the weights to fit a regular quantile regression model.
This process has to be iterated since $\beta$ in $G_0(Z_i'\beta|Z_i)$
is unknown. The iteration is initialized by using the inverse
probability estimator. For the examples in the simulation part and the
data analysis, this iteration is very fast. Convergence is achieved
usually in a few iterations for a reasonable convergence criterion.
Note that we need to estimate $G_0(Z_i'\beta|Z_i)$ at each iteration,
while no iteration is needed for Wang and Wang's algorithm.
\end{remark}


\section{Asymptotic theory}\label{sec3}
We establish the consistency and the asymptotic normality of the
estimator in this section. To derive the asymptotic properties of
the proposed estimator, we require the following regularity
assumptions. For convenience, we write the true value of $\beta$ as
$\beta_0$. 
The regularity conditions are listed as follows:
\begin{enumerate}[C1.]

\item[C1.] $T$ and $C$ are conditionally independent given the
covariate $Z$.

\item[C2.] The true value $\beta_0$ of $\beta$ is in the interior
of a bounded convex region $\mathcal{ B}$. The support $\mathcal{ Z}$ of $Z$
is bounded.\vadjust{\goodbreak}

\item[C3.] $\inf_{Z\in\mathcal{ Z}}P(Y\geq\mathcal{ T}|Z)\geq\eta_0>0$,
where $\mathcal{ T}=\mathcal{ T}_0\vee\sup_{Z\in\mathcal{ Z},\beta\in\mathcal{
B}} Z'\beta$.

\item[C4.] The conditional density functions $f_0(t|z)$ and $g_0(t|z)$
of the failure time $T$ and $C$, respectively, are uniformly
bounded away from infinity and have bounded (uniformly in $t$)
second order partial derivatives with respect to $z$.

\item[C5.] The bandwidth $h_n$ satisfies $h_n=\mathrm{O}(n^{-v})$ with
$0<v<1/2$.

\item[C6.] The kernel function $K(\cdot)\geq0$ is compactly supported
and satisfies
the Lipschitz condition of order 1, $\int K(u)=1,\int
uK(u)\,\mathrm{d}u=0,\int K^2(u)\,\mathrm{d}u<\infty$ and $\int\|u\|^2K(u)\,\mathrm{d}u<\infty$.

\item[C7.] For $\beta$ in a neighborhood of $\beta_0$,
$E[ZZ'f_0(Z'\beta|Z)]$ is positive definite.
\end{enumerate}
Assumptions C1--C4 are standard in survival analysis. Assumption C5 is needed
to ensure the consistency of the local Kaplan--Meier estimator.
Assumption C6 is
routinely made in nonparametric smoothing, and assumption C7 ensures a
unique solution for the limiting estimating equation in the
neighborhood of $\beta_0$ and is used to derive the asymptotic
properties of the estimator. Intuitively, if $\hat{G}(\beta'Z_i|Z_i)$ is
a reasonable estimator of ${G}(\beta'Z_i|Z_i)$, the consistency of
$\hat
\beta$
follows from the unbiasedness of the estimating equation~(\ref{eqyjw}). Formally,
we have the following results for the consistency of the estimator.

\begin{theorem}[(Consistency)]\label{th1} Under conditions \textup{C1--C7}, we have that
$\hat\beta_n\to\beta_0$
in probability as $n\to\infty.$
\end{theorem}

The proof of this theorem uses the uniform consistency of $\hat{G}$ as an
estimator of $G$ and is similar to that in Ying \textit{et al.}~\cite{Yingetal1995}.
Since the
criterion function is not smooth, we make use of the general theorem developed
by Chen, Linton and Van Keilegom~\cite{Chen2003} to show the asymptotic
normality of the
resulting estimator.

\begin{theorem}[(Asymptotic normality)]\label{th2} Under conditions \textup{C1--C7},
if $1/4<v<1/3$, then we have that
\[
n^{1/2}(\hat\beta_n- \beta_0)\stackrel{d}{\to} N(0,\Gamma
_1^{-1}V\Gamma_1^{-1}),
\]
where $\Gamma_1=EZZ'f_0(Z'\beta_0|Z)$ and $V=\operatorname{cov}(V_i)$ with
$V_i$ defined in Lemma \textup{\ref{lea3}} in the \hyperref[app]{Appendix}.
\end{theorem}

Note that this theorem is only for problems with a single covariate.
As in Wang and Wang~\cite{WangWang2009}, we observe that the results are not very
sensitive to $h_n$.
In practice, we can use K-folds cross validation to choose the
bandwidth. This approach works by dividing the data
set in $K$ parts, which are about equally sized. For the $k$th part,
we use the rest $K-1$ parts of the data to fit the model, and then
evaluate the quantile loss from predicting the $\tau$th conditional
quantile of $T$ on the uncensored data that are left out. Averaging
over $k=1,\ldots ,K$, we choose the $h$ that gives the minimum average
quantile loss.

The matrices $\Gamma_1$ and $V$ in
the limiting covariance matrix depend on the unknown conditional
density function
$f_0(\cdot|z)$ and $g_0(\cdot|z)$. For censored data, they may not
be estimated well nonparametrically with finite sample. Thus we
use the bootstrap resampling procedure for
inference. The validity of this procedure can be shown
following Jin \textit{et~al.}~\cite{Jinetal2003}.
We note that for the bootstrap or other resampling methods to be
feasible computationally, efficient algorithms are instrumental
because a large number of bootstrap replications are needed.
\section{Numerical study}\label{sec4}
For simulation study, we compare the estimator of Ying \textit{et al.} (YJW),
the proposed estimator (CQR),
the locally weighted censored quantile regression estimator (Lcrq)
in Wang and Wang~\cite{WangWang2009}, the estimator in Portnoy~\cite{Portnoy2003}, abbreviated
as Port,
and the estimator by Peng and Huang~\cite{PengHuang2008}, abbreviated as PH. We use
Wang and Wang's code, available on their websites, for Lcrq and
function \textit{crq} in R library \textit{quantreg} for Portnoy's and Peng and Huang's
method. YJW is implemented via the iterative method in this paper by replacing
$\hat{G}(\beta'Z_i|Z_i)$ in (\ref{eqyjw}) by $\hat{G}(\beta'Z_i)$,
which is the
Kaplan--Meier estimate for the survival function of $C$. We follow Yin
and Cai~\cite{Yin2005} and Zhou~\cite{Zhou2006} to obtain the initial value of $\beta$
by using weights $\hat{G}(Y_i)$ in (\ref{eqyjw}). A justification of
this algorithm can be found in Yin and Cai~\cite{Yin2005}.
Note that the local Kaplan--Meier estimate $\hat{G}(\beta'Z_i|Z_i)$ can
be obtained by simply modifying Wang and Wang's \textit{R} function for
their local Kaplan--Meier estimation.

We compare the mean bias (MB), the median absolute error (MAE) and
the root mean square errors (RMSE) of these procedures (Koenker~\cite{Koenker2008}).
We fix $h=0.05$
for all the simulations as suggested by Wang and Wang~\cite{WangWang2009}. Other choices
of the bandwidth were also tried.
The results are similar and are omitted to save
space. We investigate the performance at the median $\tau=0.5$
for two sample sizes $n=100$ and $200$. For Examples~\ref{ex1} and~\ref{ex2}, we have also
examined the performance at $\tau=0.7$, and the results are similar.
For each setup, the simulation is repeated 500 times. And we use 400
bootstrap replications for inference.

\begin{example}\label{ex1} We take the first example from Wang and
Wang~\cite{WangWang2009} to generate failure time
from the following i.i.d. error model
\[
T_i=b_0+b_1z_i+\varepsilon_i,\qquad i=1,\ldots ,n,
\]
where $b_0=3, b_1=5,  Z\sim U(0,1)$ and $\varepsilon_i=\eta_i-\Phi
^{-1}(\tau)$ with $\{\eta_i \}_{i=1}^n$ being i.i.d. standard normal
random variables. The censoring variable is either generated from
$U(0,14)$ or $U(0,36)$
such that about $40\%$ or $15\%$ of the observations are censored at
the median, when $\tau=0.5$ is used for generating $\varepsilon$.
\end{example}

\begin{example}\label{ex2} This example is again taken from Wang and Wang
\cite{WangWang2009}. The data are generated from
\[
T_i=b_0+b_1z_i+\bigl(0.2+a(z_i-0.5)^2\bigr)\varepsilon_i,\qquad i=1,\ldots
,n,
\]
where $b_0=2$, $b_1=1$, $z_i \sim N(0,1)$ and
$\varepsilon_i=\eta_i-\Phi^{-1}(\tau)$ with $\{\eta_i \}_{i=1}^n$
being i.i.d. standard normal random variables. Here $a$ takes the value
$0$, $0.5$
and $2$ to indicate no, median and strong deviation from the global linearity
assumption. The
censoring variable $C_i$ is generated from $U(0,7)$ or $U(0,18)$ to
give $40\%$
and $15\%$ censoring at the median for $a=2$.
\end{example}

\begin{example}\label{ex3} The model to generate $T_i$ is the same as
Example~\ref{ex2} with $b_0=1$ and $a=2$.
The censoring variable $C_i$ is generated from a mixture
of distributions. Specifically, if $z_i<1$, $C_i$ is generated from
$U(0,4)$; otherwise, $C_i$ is generated from $U(0,8)$. This scheme
gives about
$30\%$ censoring at the median.
\end{example}

\begin{example}\label{ex4} This model is similar to Example~\ref{ex2} with
$a=2$. However, the censoring time $C_i$ is generated from the
following model
$C_i=A+bz_i+\eta_i, i=1,\ldots ,n,$
where $b=0, 0.5, 1$ indicates a different level of dependence of the censoring
time on the covariates. The random variable $\eta_i$ follows the
standard normal distribution, and $A$ is either 1.35 or 2.6, such that when
$b=1$, about $40\%$ or $15\%$ observations are censored.
\end{example}

For Examples~\ref{ex1} and~\ref{ex2}, the censoring time and the survival time are independent.
In Example~\ref{ex1}, all the conditional
quantiles are linear functions of the covariates. Example~\ref{ex2} gives a model
with only the $\tau$th quantile being a linear function when $a\not= 0$.
Note that Wang and Wang~\cite{WangWang2009} have shown that Portnoy's approach gives
biased estimates for the coefficients when $a=2$. For Examples~\ref{ex3} and~\ref{ex4},
$T$ and $C$ are not independent, but they
are conditionally independent given the covariates. Example~\ref{ex3} uses a
mixture distribution to generate censoring time, while
in Example~\ref{ex4}, a linear dependence of $C$ on the covariates is used. Two
different censoring rates are examined for Examples~\ref{ex1},~\ref{ex2} and~\ref{ex4}.

It is not difficult to see that the initial estimate $\beta^{(0)}$ is
also $\sqrt{n}$-consistent. However, we observe empirically that
it is less efficient than the final
estimate after iteration. For example,
in Example~\ref{ex2} when $n=100$, $a=2$ and censoring rate is $40\%$,
the RMSEs of $\beta^{(0)}$ are 0.438 and
0.697 when censoring rate is $40\%$, and 0.246 and 0.450 when about
$15\%$
of the data are censored. A related comparison was made in Yin, Zeng
and Li~\cite{Yingetal2008}.
The results for the other methods
are summarized in Tables~\ref{tab1} and~\ref{tab11} when $\tau=0.5$ and $n=100$.
The results for $\tau=0.7$ or $n=200$ are qualitatively similar and
thus are omitted. We have the following observations. First, CQR
outperforms YJW in general, especially when the unconditional
independence is violated.
Second, when the global linearity holds, Port and PH generally
outperform CQR and Lcrq, although by a small margin. When the global
linearity is mildly violated, PH and Port both perform competitively
with CQR and Lcrq. This demonstrates the robustness of PH and Port.
However, when this assumption is severely violated, CQR and Lcrq
perform better in general, especially in terms of RMSE. However, how
much improvement can be expected is likely dependent on a number of
factors, such as the censoring mechanism and rates. Third, CQR and Lcrq
have similar performance across all the simulations. The difference
between these two approaches is usually negligible. Fourth, when the
censoring is low ($15\%$), CQR performs competitively compared to Lcrq,
indicating its robustness with respect to the required sample size for
estimating the local Kaplan--Meier curve. We conclude that when the
global linearity is violated, and $C$ is not unconditionally
independent of $T$, the
proposed method is preferred over YJW, Port and PH.

%
\begin{table}
\caption{Simulation results for Examples \protect\ref{ex1} and \protect\ref{ex2}} \label{tab1}
\begin{tabular*}{\textwidth}{@{\extracolsep{\fill}}lllld{2.3}d{2.3}cccc@{}}
\hline
&&&&\multicolumn{2}{l}{Bias}&\multicolumn{2}{l}{MAE}&\multicolumn{2}{l@{}}{RMSE}\\[-6pt]
&&&&\multicolumn{2}{l}{\hrulefill}&\multicolumn{2}{l}{\hrulefill}&\multicolumn{2}{l@{}}{\hrulefill}\\
Ex.& \multicolumn{1}{l}{$c\%$}& \multicolumn{1}{l}{$a$}&&
\multicolumn{1}{l}{$\beta_0$}&
\multicolumn{1}{l}{$\beta_1$}&
\multicolumn{1}{l}{$\beta_0$}&
\multicolumn{1}{l}{$\beta_1$}&
\multicolumn{1}{l}{$\beta_0$}&
\multicolumn{1}{l}{$\beta_1$}\\
\hline
\ref{ex1}
&$40\%$&&YJW&0.010 & -0.015& 0.228 & 0.431 & 0.341 & 0.644\\
& & &CQR& -0.007 & -0.092 & 0.211 & 0.392 & 0.305 & 0.583\\
&& &Lcqr& -0.009 & -0.018 & 0.196 & 0.375 &0.299 & 0.554\\
&& &Port& -0.042 & -0.003 &0.198 & 0.381 & 0.299 & 0.548\\
&& &PH& 0.019 & 0.005 & 0.203 & 0.380 & 0.299 & 0.558\\[3pt]
&$15\%$&&YJW& -0.016 & 0.023 & 0.191 & 0.330 & 0.293 & 0.508\\
& & &CQR& -0.012 & -0.012 & 0.186 & 0.317 & 0.280 & 0.481\\
&& &Lcqr& -0.013 & 0.005 & 0.186 & 0.311 & 0.281 & 0.485 \\
&& &Port& -0.059 & 0.007 & 0.179 & 0.298 & 0.279 & 0.469\\
&& &PH& 0.001 & 0.007 & 0.186 & 0.305 & 0.277 & 0.477\\[6pt]
\ref{ex2}&$40\%$&2
&YJW&0.007 & -0.009 & 0.159 & 0.326& 0.249 & 0.537\\
&& &CQR& -0.060 & -0.030 & 0.139 & 0.267 & 0.211 & 0.393 \\
&& &Lcqr& -0.053 & 0.008 & 0.144 & 0.272 & 0.215 & 0.406 \\
&& &Port& -0.021 & -0.010 & 0.164 & 0.308 & 0.224 & 0.443 \\
&& &PH& 0.058 & -0.119  &0.166 & 0.304 & 0.235 & 0.460\\[3pt]
&&0.5&YJW& -0.001 & 0.012 & 0.058 & 0.125 & 0.089 & 0.188 \\
& & &CQR& -0.044 & -0.016 & 0.064 & 0.106 & 0.095 & 0.154\\
&& &Lcqr& -0.021 & 0.009 & 0.058 & 0.106 & 0.088 & 0.164\\
&& &Port& -0.020 & 0.012& 0.058 & 0.109 & 0.089 & 0.168\\
&& &PH& 0.013 & -0.014 & 0.057 & 0.106 & 0.087 & 0.169\\[3pt]
&&0&YJW&0.003 & 0.008 & 0.022 & 0.025 & 0.033 & 0.040\\
& & &CQR& -0.022 & -0.012 & 0.028 & 0.024 & 0.041 & 0.037\\
&& &Lcqr& -0.004 & -0.002 & 0.021 & 0.022 & 0.031 & 0.034\\
&& &Port& -0.010 & 0.001 & 0.021 & 0.021& 0.032 & 0.033 \\
&& &PH& 0.003 & 0.001 & 0.021 & 0.021 & 0.031 & 0.033\\[3pt]
&$15\%$&2 &YJW& -0.006 & 0.012 & 0.146 & 0.292 & 0.212 & 0.425\\
&& &CQR& -0.024 & -0.004 & 0.134 & 0.267 & 0.202 & 0.393 \\
&& &Lcqr& -0.023 & 0.013 & 0.138 & 0.271 & 0.202 & 0.396\\
&& &Port& -0.051 & 0.039 & 0.141 & 0.277 & 0.214 & 0.410\\
&& &PH& 0.021 & -0.038 & 0.143 & 0.288 & 0.208 & 0.406\\[3pt]
&&0.5&YJW& 0.001 & 0.000 & 0.055 & 0.100 & 0.082 & 0.150 \\
& & &CQR& -0.011 & -0.008 & 0.055 & 0.092 & 0.082 & 0.141\\
&& &Lcqr& -0.005 & 0.000 & 0.052 & 0.092 & 0.081 & 0.144\\
&& &Port& -0.025 & 0.014& 0.056 & 0.091 & 0.085 & 0.144\\
&& &PH& 0.008 & -0.008 & 0.054 & 0.094 & 0.081 & 0.144\\[3pt]
&&0&YJW&0.001 & 0.003 & 0.018 & 0.019 & 0.027 & 0.029 \\
& & &CQR& -0.007 & -0.003 & 0.020 & 0.019 & 0.028 & 0.028\\
&& &Lcqr& -0.001 & -0.000 & 0.019 & 0.018 & 0.027 & 0.027\\
&& &Port& -0.011 & -0.000 & 0.020 & 0.018 & 0.029 & 0.027\\
&& &PH& 0.001 & 0.000 & 0.018 & 0.018 & 0.026 & 0.027\\
\hline
\end{tabular*}
\end{table}

%
\begin{table}
\caption{Simulation results for Examples \protect\ref{ex3} and \protect\ref{ex4}}
\label{tab11}
\begin{tabular*}{\textwidth}{@{\extracolsep{\fill}}lllld{2.3}d{2.3}cccc@{}}
\hline
&&&&\multicolumn{2}{l}{Bias}&\multicolumn{2}{l}{MAE}&\multicolumn{2}{l@{}}{RMSE}\\[-6pt]
&&&&\multicolumn{2}{l}{\hrulefill}&\multicolumn{2}{l}{\hrulefill}&\multicolumn{2}{l@{}}{\hrulefill}\\
Ex.& \multicolumn{1}{l}{$c\%$}& \multicolumn{1}{l}{$a$}&&
\multicolumn{1}{l}{$\beta_0$}&
\multicolumn{1}{l}{$\beta_1$}&
\multicolumn{1}{l}{$\beta_0$}&
\multicolumn{1}{l}{$\beta_1$}&
\multicolumn{1}{l}{$\beta_0$}&
\multicolumn{1}{l}{$\beta_1$}\\
\hline
\ref{ex3}&$30\%$&&YJW&0.035 & 0.544 & 0.124 & 0.482 & 0.206 & 0.810\\
&&&CQR & 0.005 & -0.023 & 0.115 & 0.223 & 0.164 & 0.325\\
&&&Lcrq & 0.011 & 0.005 & 0.109 & 0.234 & 0.164 & 0.333\\
&&&Port & 0.020 & 0.017 & 0.124 & 0.257 & 0.190 & 0.367 \\
&&&PH & 0.088 & -0.058 & 0.128 & 0.261 & 0.229 & 0.388\\[5pt]
\ref{ex4}&$30\%$&1 &YJW& -0.115 & 0.421 & 0.121 & 0.411 & 0.159 & 0.473\\
&& &CQR& -0.077 & 0.081 & 0.095 & 0.137 & 0.137 & 0.217\\
&& &Lcqr& -0.080 & 0.083 & 0.095 & 0.138 & 0.139 & 0.218 \\
&& &Port& -0.029 & 0.024 & 0.086 & 0.140 & 0.131 & 0.226\\
&& &PH&0.017 & -0.020 & 0.085 & 0.148 & 0.130 & 0.222\\[3pt]
&&0.5&YJW& -0.063 & 0.227 & 0.0944 & 0.228 & 0.137 & 0.331\\
& & &CQR& -0.055 & -0.011 & 0.090 & 0.155 & 0.131 & 0.217\\
&& &Lcqr& -0.051 & -0.010 & 0.090 & 0.153 & 0.129 & 0.222 \\
&& &Port& -0.018 & -0.022 & 0.093 & 0.178 & 0.136 & 0.251\\
&& &PH& 0.035 & -0.089 & 0.090 & 0.189 & 0.141 & 0.267\\[3pt]
&&0&YJW&0.000 & -0.020 & 0.096 & 0.207 & 0.143 & 0.311\\
& & &CQR& -0.042 & -0.081 & 0.088 & 0.191 & 0.127 & 0.264\\
&& &Lcqr&-0.040 & -0.056 & 0.086 & 0.199 & 0.129 & 0.279 \\
&& &Port& -0.020 & -0.037 & 0.088 & 0.207 & 0.135 & 0.306\\
&& &PH& 0.030 & -0.111 & 0.092 & 0.234 & 0.135 & 0.319\\[3pt]
&$15\%$&1&YJW&-0.020 & 0.226 & 0.079 & 0.226 & 0.117 & 0.323\\
& & &CQR& -0.014 & 0.023 & 0.078 & 0.148 & 0.114 & 0.215 \\
&& &Lcqr& -0.015 & 0.021 & 0.076 & 0.146 & 0.114 & 0.213 \\
&& &Port& -0.038 & 0.029 & 0.084 & 0.155 & 0.125 & 0.227 \\
&& &PH& 0.007 & -0.007 & 0.078 & 0.152 & 0.118 & 0.223\\[3pt]
&&0.5&YJW&-0.017 & 0.157 & 0.081 & 0.187 & 0.123 & 0.290\\
&& &CQR& -0.004 & -0.015 & 0.078 & 0.141 & 0.119 & 0.218 \\
&& &Lcqr&-0.002 & -0.021 & 0.077 & 0.144 & 0.119 & 0.216 \\
&& &Port& -0.034 & 0.008 & 0.088 & 0.145 & 0.130 & 0.225\\
&& &PH& 0.014 & -0.036 & 0.081 & 0.143 & 0.124 & 0.226\\[3pt]
&&0&YJW&0.000 & -0.001 & 0.083 & 0.167 & 0.123 & 0.260\\
& & &CQR&-0.000 & -0.052 & 0.080 & 0.156 & 0.119 & 0.235 \\
&& &Lcqr& 0.000 & -0.058 & 0.079 & 0.160 & 0.119 & 0.237\\
&& &Port& -0.031 & -0.021 & 0.085 & 0.157 & 0.126 & 0.237\\
&& &PH& 0.017 & -0.077 & 0.079 & 0.166 & 0.123 & 0.248\\
\hline
\end{tabular*}
\end{table}

We assess the performance of the bootstrap inference procedure by
comparing it
to the bootstrap percentile inference
procedure developed in Wang and Wang~\cite{WangWang2009}. For brevity, we only report
the result for Example~\ref{ex2} when $a=2$, and the censoring rate is $40\%$,
and for
Example~\ref{ex3}. We record the
empirical coverage probability (ECP) and the empirical mean length
(EML) of the resulting confidence intervals in Table~\ref{tab3}.
The nominal level used is
0.95. For these two examples, both CQR and Lcrq give coverage probabilities
close to the
nominal level with comparable average empirical lengths.

%
\begin{table}
\caption{Performance of the inference procedure.
For Example \protect\ref{ex2}, only the result for $a=2$ and $40\%$ censoring is
presented} \label{tab3}
\begin{tabular*}{\textwidth}{@{\extracolsep{\fill}}lllllllllll@{}}
\hline
&&&\multicolumn{4}{l}{CQR}&\multicolumn{4}{l}{Lcrq}\\[-6pt]
&&&\multicolumn{4}{l}{\hrulefill}&\multicolumn{4}{l@{}}{\hrulefill}\\
&&&\multicolumn{2}{l}{ECP}&\multicolumn{2}{l}{EML}&\multicolumn{2}{l}{ECP}&\multicolumn{2}{l@{}}{EML}\\
[-6pt]
&&&\multicolumn{2}{c}{\hrulefill}&\multicolumn{2}{c}{\hrulefill}&\multicolumn{2}{c}{\hrulefill}&\multicolumn{2}{c@{}}{\hrulefill}\\
\multicolumn{1}{@{}l}{Ex.}&$n$&$\tau$&$\beta_0$&$\beta_1$&$\beta_0$&$\beta_1$&$\beta
_0$&$\beta_1$&$\beta_0$&$\beta_1$\\
\hline
2&100&0.5&0.948 & 0.930 & 0.800 & 1.541 &0.938 & 0.950 & 0.806 &1.564\\
&&0.7& 0.938 & 0.944 & 0.862 & 1.678 & 0.956 & 0.954 & 0.883 & 1.740\\
&200&0.5& 0.928 & 0.934 & 0.533 & 1.031 & 0.946 & 0.938 & 0.586 &
1.147\\
&&0.7&0.930 & 0.932 & 0.567 & 1.102 & 0.934 & 0.948 & 0.596 &1.174
\\[5pt]
3&100&0.5& 0.910 & 0.930 & 0.612 & 1.195& 0.942 & 0.924 & 0.596 &
1.196\\
&&0.7&0.916 & 0.928 & 1.037 & 2.086 & 0.928 & 0.936 & 0.956 & 1.966\\
&200&0.5 & 0.916 & 0.916 & 0.423 & 0.839 & 0.954 & 0.950 & 0.424 &
0.872\\
&&0.7 &0.934 & 0.928 & 0.772 & 1.490 & 0.946 & 0.944 & 0.715 & 1.445
\\
\hline
\end{tabular*}
\end{table}

\begin{table}
\caption{Multi-dimensional covariates when $n=200$. The standard
errors (SE)
are reported in parentheses. Note that the SEs of MAE and RMSE are
computed via bootstrapping 1000 replications} \label{tab14}
\begin{tabular*}{\textwidth}{@{\extracolsep{\fill}}lld{2.4}d{2.4}d{2.4}d{2.4}d{2.4}d{2.4}@{}}
\hline
&&\multicolumn{2}{l}{Bias}&\multicolumn{2}{l}{MAE}&\multicolumn{2}{l@{}}{RMSE}\\[-6pt]
&&\multicolumn{2}{l}{\hrulefill}&\multicolumn{2}{l}{\hrulefill}&\multicolumn{2}{l@{}}{\hrulefill}\\
\multicolumn{1}{@{}l}{$d$}&\multicolumn{1}{l}{$c\%$} &
\multicolumn{1}{l}{$\beta_0$}&
\multicolumn{1}{l}{$\beta_1$}&
\multicolumn{1}{l}{$\beta_0$}&
\multicolumn{1}{l}{$\beta_1$}&
\multicolumn{1}{l}{$\beta_0$}&
\multicolumn{1}{l@{}}{$\beta_1$}\\
\hline
1&$40\%$&-0.010 & -0.034 & 0.152 & 0.277 & 0.220 & 0.413 \\
& & (0.220) & (0.412) & (0.006) & (0.009) & (0.005) & (0.009)\\
&$15\%$&0.001 & -0.010 & 0.127 & 0.227 & 0.191 & 0.342 \\
& & (0.191) & (0.342) & (0.006) & (0.010) & (0.005) & (0.008)\\[3pt]
2&$40\%$&-0.070 & -0.260 & 0.199 & 0.338 & 0.294 & 0.493 \\
& & (0.286) & (0.419) & (0.009) & (0.012) & (0.007) & (0.010)\\
&$15\%$& -0.027 & -0.063 & 0.177 & 0.245 & 0.258 & 0.350 \\
& & (0.256) & (0.344) & (0.006) & (0.007) & (0.006) & (0.008)\\[3pt]
3&$40\%$&-0.100 & -0.339 & 0.247 & 0.378 & 0.369 & 0.535 \\
& & (0.355) & (0.414) & (0.009) & (0.010) & (0.008) & (0.011)\\
&$15\%$&-0.032 & -0.096 & 0.211 & 0.236 & 0.315 & 0.353 \\
& & (0.314) & (0.341) & (0.008) & (0.007) & (0.007) & (0.007)\\[3pt]
4&$40\%$&-0.087 & -0.414 & 0.285 & 0.452 & 0.438 & 0.602 \\
& & (0.430) & (0.438) & (0.011) & (0.016) & (0.011) & (0.013)\\
&$15\%$&-0.029 & -0.130 & 0.238 & 0.229 & 0.355 & 0.359 \\
& & (0.354) & (0.334) & (0.011) & (0.009) & (0.007) & (0.010)\\[3pt]
5&$40\%$&-0.089 & -0.441 & 0.364 & 0.465 & 0.493 & 0.616 \\
& & (0.485) & (0.430) & (0.013) & (0.019) & (0.011) & (0.014)\\
&$15\%$& -0.048 & -0.125 & 0.254 & 0.245 & 0.387 & 0.367 \\
& & (0.384) & (0.345) & (0.009) & (0.009) & (0.008) & (0.009)
\\\hline
\end{tabular*}
\end{table}

Since CQR relies on the local Kaplan--Meier estimate, it is of great
interest to
see how it performs when $Z$ is multi-dimensional. To this end, we use
the same model in Example~\ref{ex1}, but add independent standard uniform
covariates $z^{(2)},\ldots ,z^{(d)}$ to $z$. Thus, the coefficients
associated with these additional covariates are zero. We use bandwidth
0.05, 0.1, 0.2, 0.3, 0.4, respectively, when $d=1, 2,\ldots ,5$. In Table
\ref{tab14}, we see that with $n=200$, CQR seems to give unbiased
estimates when $d=1, 2$ and $3$ with low censoring $15\%$. However, it
can only be applied up to $d=2$ if censoring rate is as high as $40\%$.

%

\section{Data analysis}\label{sec5}

\begin{sidewaysfigure}

\includegraphics{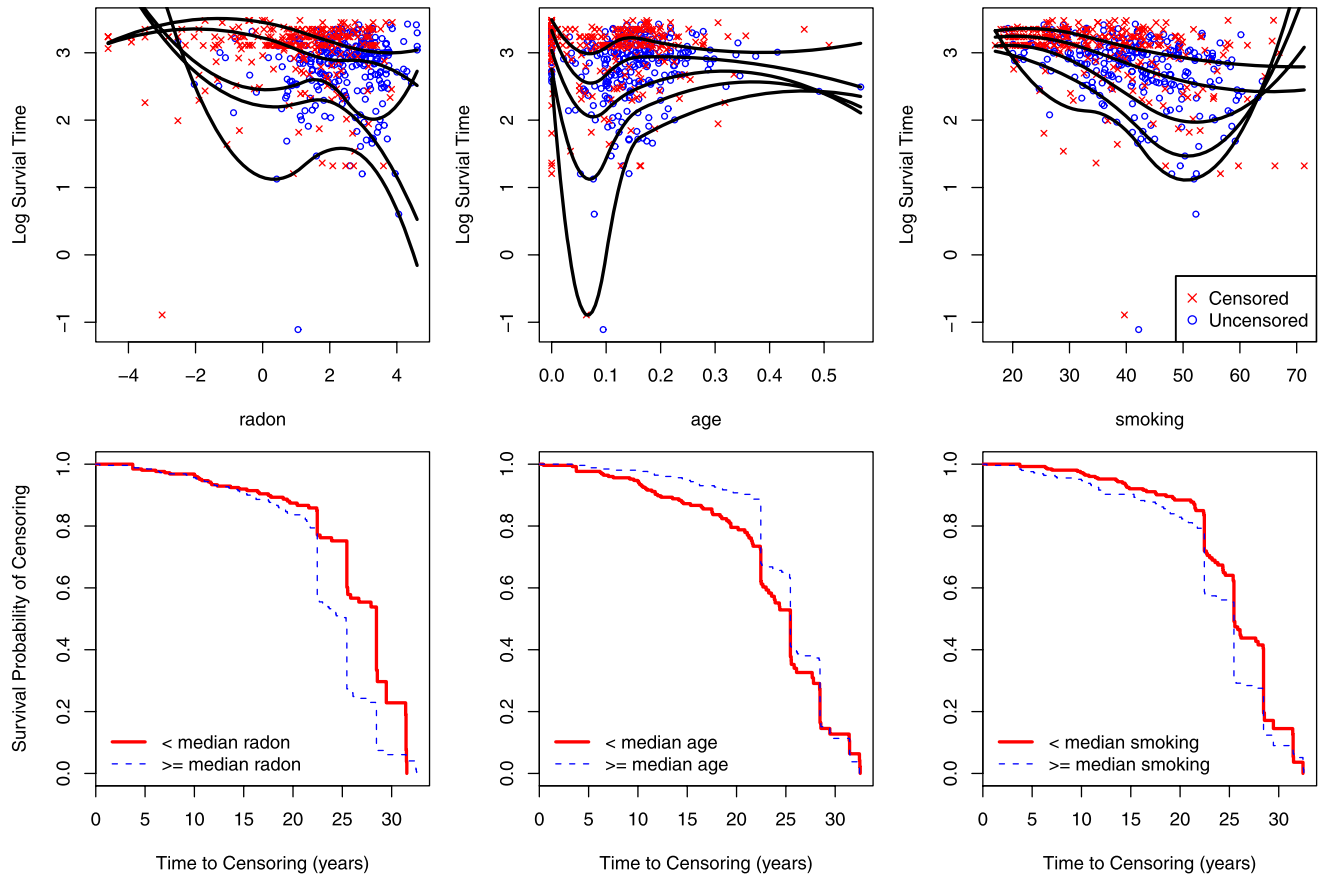}

\caption{Colorado miners cohort data. Top row: The scatter plots of
the log survival time versus the covariates. The marginally fitted log
survival times against each covariate at quantiles 0.01, 0.05, 0.1,
0.3, 0.5 (the solid lines from the bottom to the top) are also plotted.
Bottom row: The fitted Kaplan--Meier survival curves for the censoring
time when covariates are dichotomized.}
\label{fig1}
\end{sidewaysfigure}

As an example, we apply the proposed method to the Colorado
Plateau uranium miners cohort data (Lubin \textit{et al.}~\cite{L1994}, Langholz and
Goldstein~\cite{Lang1996}). The major interest of this study is to assess the
effect of smoking on the rate of median lung cancer. This data set
consists of 3347 Caucasian male miners who worked underground for at
least one month in the uranium mines of the Colorado Plateau area.
In total, there are 258 miners who died of lung cancer. Apart from
the failure time, information of the age, the cumulative radon
exposure and cumulative smoking in number of packs is available. In
our study, we randomly choose 258 miners who are censored and all
the miners who experience the lung cancer. We use this scheme to yield a
median censoring scenario, suggested by the simulation studies. This data
analysis means to illustrate the difference between different approaches.
The scatter plots of the log
survival time are presented at the top row of Figure~\ref{fig1}.
Let $Z_1$ be the
logarithm of the cumulative radon exposure in 100 working level
months, $Z_2$ be the cumulative smoking in 1000 packs and $Z_3$ be
the age at entry to the study. To explore the dependence of the log survival
time against these covariates, we fit three separate marginal models
using polynomial B-splines to approximate the effects of radon, age and
smoking, respectively. In Figure~\ref{fig1}, we plot the estimated log survival
time against the three covariates at quantiles $\tau
=0.01, 0.05, 0.1, 0.3, 0.5$. Strong non-linearity is present,
especially for lower quantiles.
When $\tau=0.5$, the log survival time is approximately linear. These facts
suggest that the global linearity assumption may not hold.
To further examine whether unconditional
independence of the survival time and censoring time is appropriate,
we fit the Cox model to the censoring time with respect to the
covariates. The two covariates radon and age are both found
significant from zero with p-values less than $10^{-3}$. This
indicates that the unconditional independent assumption needed for
YJW may be inappropriate for this data. Graphically, the
dependence of the censoring time on the covariates can be seen from
Figure~\ref{fig1}, where Kaplan--Meier estimates of the survival
functions, dichotomized by the median of these covariates, are
plotted. From the figure, an observation is more likely to be
censored at an earlier time if radon
%
is high, age is young or the
subject smoked less. Formal log-rank tests by dichotomizing the
covariates also indicate that radon ($p$-value $<10^{-3}$) and smoking
($p$-value $0.03$) are highly correlated with the censoring time,
while age ($p$-value $0.08$) is not significant. Note that these log-rank
tests only investigate these covariates marginally.

%
\begin{table}
\caption{The fitted coefficients of the censored quantile regression at
the median and the $95\%$
confidence intervals for the Colorado Plateau uranium miners cohort data}\label{tab7}
\begin{tabular*}{\textwidth}{@{\extracolsep{\fill}}llll@{}}
\hline
& CQR &YJW&Lcrq
\\
\hline
Radon ($\times10^{-2}$)&$-13.05_{(-28.85,-9.84)}$&$-2.95_{(-12.4,3.23)}$&
$-25.60_{(-34.45,-13.55)}$
\\
Smoking ($\times10^{-3}$)&\phantom{0}$-0.38_{(-3.73,3.30)}$&\phantom{0.}$0.04_{(-5.19,3.58)}$&
\phantom{00.}$0.65_{(-4.91,2.86)}$
\\
Age&\phantom{0}$-1.60_{(-2.27,-1.28)}$&$-2.01_{(-2.61,-1.44)}$&
\phantom{0}$-2.12_{(-2.76,-1.39)}$
\\
\hline
\end{tabular*}
\end{table}

Since we have three continuous covariates, we use the three-dimensional
kernel after standardizing the covariates, which is the product of
three bi-quadratic kernels for radon, smoking and age. We investigate
the median log survival time on the three covariates. For Lcrq and CQR,
we use the same bandwidth for the three univariate kernels and apply
10-fold cross validation to choose the optimal bandwidth. The results
are summarized in Table~\ref{tab7}. It is seen that Ying \textit{et al.}
estimate age as the only significant variable, while CQR and Lcrq both
estimate age and radon as significant variables. The result of Ying
\textit{et al.}'s approach in this example is problematic due to the
dependence of the censoring time on the covariates. The $95\%$
confidence intervals are obtained by using the bootstrap percentile
approach (Wang and Wang~\cite{WangWang2009}) using 1000 bootstrap repetitions. The
fact that we only use a random sample for the censored data suggests
that the results, especially the numerical ones, should be interpreted
with certain caution. 

\section{Conclusion}\label{sec6}
We have proposed a novel extension of Ying, Jung and Wei's median
regression to quantile regression. Our model is more flexible in
that only conditional independence of the survival time and
censoring time are assumed. Moreover, we have proposed a new and fast
fitting algorithm, applicable to Ying \textit{et al.}'s median regression
model, making use of the efficient quantile regression code
developed by Koenker. Therefore, resampling based inference
procedure can be efficiently implemented. We have compared our
estimator to the approaches developed in Portnoy~\cite{Portnoy2003}, Peng and
Huang~\cite{PengHuang2008} and Wang and Wang~\cite{WangWang2009}. The simulation results show
that the new method is useful and may have certain advantages over
the other methods, especially when the global linearity is violated or the
unconditional independence of $C$ and $T$ does not hold.

Identifiability remains a serious issue in censored quantile
regression, particularly so when~$\tau$ is close to 1 or 0 (Peng and
Huang~\cite{PengHuang2008}, Wang and Wang~\cite{WangWang2009}). In practice, we recommend to
choose~$\tau$ in the inference range of interest. Another limitation of the
current method is the requirement of estimating $G_0(\cdot|Z)$, which
inevitably suffers from the curse of dimensionality if~$Z$ is
multi-dimensional. In this case, it may be more attractive to handle
$G_0(\cdot|Z)$ by using, for example, the Cox model or the single-index
model. This line of research merits further investigation. Furthermore,
the Kaplan--Meier estimates, even for a global one, may be unstable at
the right tails. The technique in Zhou~\cite{Zhou2006} may be used to improve
the stability of these estimates.


\begin{appendix}\label{app}
\section*{Appendix}
For convenience, we write $\|\beta\|$ as the Euclidean norm of a
finite dimensional vector $\beta$ and $\|G(\cdot)\|_{\infty}$ as the
supreme of the absolute value of a function $G(\cdot)$.
First, we cite Theorem 2.1 in Gonzalez and Cadarso~\cite{G1994}.

\renewcommand{\thelemma}{A.\arabic{lemma}}
\begin{lemma}\label{lea1} Assume that conditions \textup{C4--C6} hold, then
\[
\|\hat{G}-G_0\|_{\infty}=\sup_t\sup_z\bigl|\hat
{G}(t|z)-G_0(t|z)\bigr|=\mathrm{O}_p\bigl((\log
n)^{1/2}n^{-1/2+v/2}+n^{-2v}\bigr).
\]
\end{lemma}

\subsection{\texorpdfstring{Proof of Theorem \protect\ref{th1}}{Proof of Theorem 1}}
 Let $\widetilde{M}_n(\beta)=\sum_{i=1}^n(\tau-F_0(Z_i'\beta|Z_i))Z_i.
$ It follows from the similar arguments as in Ying \textit{et al.}~\cite{Yingetal1995} and Lemma~\ref{lea1} that
\setcounter{equation}{0}
\begin{equation}\label{A1}
\sup_{\beta\in\mathcal{ B}}
n^{-1}|M_n(\beta)-\widetilde{M}_n(\beta)|=\mathrm{o}(1) \qquad \mbox{a.s.}
\end{equation}

From assumption C7, $A_n(\beta)=\frac1n\frac{\partial
\widetilde{M}_n(\beta) }{\partial\beta}=-\frac1n\sum_{i=1}^n
Z_iZ_i'f_0(Z_i'\beta|Z_i)$ is negative definite with probability
one for $\beta$ in a small neighborhood of $\beta_0$. In addition,
$\widetilde{M}_n(\beta_0)=0$. Therefore,
$n^{-1}\widetilde{M}_n(\beta)$ is bounded away from zero. This
argument, together with (\ref{A1}), yields that $\hat\beta_n\to\beta_0$
in probability as $n\to\infty.$

\subsection{\texorpdfstring{Proof of Theorem \protect\ref{th2}}{Proof of Theorem 2}}
To prove
Theorem~\ref{th2}, we exploit Theorem 2 in Chen \textit{et al.}~\cite{Chen2003} by verifying
their conditions (2.1)--(2.4), (2.5$'$) and (2.6$'$). For convenience,
write $M_n(\beta,G)=\frac1n\sum_{i=1}^n m_i(\beta,G),$ where
$m_i(\beta,G)= Z_i\{\frac{I(Y_i\geq
Z_i'\beta)}{G(Z_i'\beta|Z_i)}-(1-\tau)\} $ and the function class $
\mathcal{ G}$ that involves the true $G_0$ as the set of $G$, such that
$G$ has a density function $g$, $\inf_{z\in\mathcal{ Z}} G(\mathcal{
T}|z)\geq\eta_0$ and $g(\cdot|z)$ is bounded away from infinity
uniformly in $t$ and $z\in\mathcal{ Z}$. Then
$M(\beta,G)=Em_i(\beta,G)=EZ_i\{
\frac{(1-F_0(Z_i'\beta|Z_i))G_0(Z_i'\beta|Z_i)}{G(Z_i'\beta
|Z_i)}-(1-\tau)\},$
where the expectation operator is taken with respect to the marginal
distribution function of $Z_i$ and thus
$M(\beta_0,G_0)=0.$

\begin{lemma}\label{lea2} For any positive value $\xi_n=\mathrm{o}(1)$, we have that
\[
\sup_{\|\beta-\beta_0\|\leq\xi_n,\|G-G_0\|_\infty\leq
\xi_n}\|M_n(\beta,G)-M(\beta,G)-M_n(\beta_0,G_0)\|=\mathrm{o}_p(n^{-1/2}).
\]
\end{lemma}

\begin{pf}
Let 
$\eta_1=\sup_{z\in\mathcal{ Z}}\|z\|^2\vee1$ and $\eta_2=\sup_{G\in
\mathcal{ G},z\in\mathcal{ Z},t\leq\mathcal{ T}}(f_0(t|z)+g(t|z))<\infty$ from
assumption C4. For any $(\beta,G)\in\mathcal{ B}\times\mathcal{ G}$ and
$(\beta^*,G^*)\in\mathcal{ B}\times\mathcal{ G}$, we have that
$\|m(\beta,G)-m(\beta^*,G^*)\|^2\leq2(U_1+U_2+U_3)$ where
\begin{eqnarray}\label{A2}
U_1&=&\bigl\|ZG(Z'\beta|Z)^{-1}\bigl(I(Y\geq Z'\beta)-I(Y\geq
Z'\beta^*)\bigr)\bigr\|^2\nonumber\\
&\leq&\eta_3|I(Y\geq
Z'\beta)-I(Y\geq Z'\beta^*)|,
\nonumber
\\[-8pt]
\\[-8pt]
\nonumber
U_2&=&\bigl\|ZI(Y\geq
Z'\beta^*)\bigl(G(Z'\beta|Z)^{-1}-G^*(Z'\beta|Z)^{-1}\bigr)\bigr\|^2\leq
\eta_4\|G-G^*\|_{\infty}^2,\\
U_3&=&\bigl\|ZI(Y\geq
Z'\beta^*)\bigl(G^*(Z'\beta|Z)^{-1}-G^*(Z'\beta^*|Z)^{-1}\bigr)\bigr\|^2\leq
\eta_5\|\beta-\beta^*\|^2,\nonumber
\end{eqnarray}
where $\eta_3,\eta_4,\eta_5$ are some positive constants, only
depending on $\eta_k$ $(k=0,1,2)$. It follows from~(\ref{A2}) that
$E(\sup_{\|\beta-\beta^*\|\leq\xi_n}U_1) \leq\eta_1
\eta_2\eta_3\xi_n $ and that
\begin{equation}\label{A3}
\sup_{\|\beta-\beta_0\|\leq\xi_n,\|G-G_0\|_\infty\leq
\xi_n}\|M(\beta,G)-M(\beta^*,G^*)\|^2\leq\eta_6\xi_n
\end{equation}
for some constant $\eta_6\geq0$ as $n$ is sufficiently large.

Therefore, condition (3.2) of Chen \textit{et al.}~\cite{Chen2003} holds with $r=2$
and $s_j=1/2$. Similarly to the arguments used in (\ref{A2}), condition
(3.1) in Chen \textit{et al.}~\cite{Chen2003} can be also verified. Now we verify
their condition (3.3). Let $N(\eta,\mathcal{ G},\|\cdot\|_{\infty})$ be
the covering numbers (van der Vaart and Wellner~\cite{VWellner1996}, page 83) for
the function class $\mathcal{ G}$ under the metrics
$\|\cdot\|_{\infty}$. An application of Theorem 2.7.1 in van der
Vaart and Wellner~\cite{VWellner1996} from assumptions C4 and C2 gives that
the logarithm of
the covering number of $\mathcal{ G}$ is bounded by $ K\eta^{-1/2}$ for
$\eta\leq1$, where $K$ is some constant, not depending on $n$. When
$\eta\geq1$, it follows from the definition of covering numbers
that $\log N (\eta,\mathcal{ G},\|\cdot\|_{\infty})=0$, which yields
that
\[
\int_0^\infty\{\log N (\eta^2,\mathcal{
G},\|\cdot\|_{\infty})\}^{1/2}\,\mathrm{d}\eta\le
\int_0^1K^{1/2}\eta^{-1/2}\,\mathrm{d}\eta<\infty.
\]
It then follows easily from Theorem 3 of Chen \textit{et al.}~\cite{Chen2003} that
Lemma~\ref{lea2} holds.
\end{pf}



To apply Theorem 2 in Chen \textit{et al.}~\cite{Chen2003}, we define
$\Gamma_1(\beta_0,G_0)$ as the first derivative function of $
M(\beta,G_0)$ with respect to $\beta$ evaluated at $\beta=\beta_0$.
For all $\beta\in\mathcal{ B}$, we define the functional derivative of
$M(\beta,G)$ at $G_0$ in the direction $[G-G_0]$ as
\[
\Gamma_2(\beta,G_0)[G-G_0]=\lim_{\eta\to
0}\frac1\eta\bigl[M\bigl(\beta,G_0+\eta(G-G_0)\bigr)-M(\beta,G_0)\bigr].
\]

\begin{lemma}\label{lea3} Assume that the conditions in Theorem~\ref{th2} hold, then
\[
n^{1/2}\bigl(M_n(\beta_0,G_0)+\Gamma_2(\beta_0,G_0)[\hat
{G}-G_0]\bigr)\stackrel
{d}{\to}N(0,V),
\]
where $V=\operatorname{cov}(V_i)$ with
$V_i= m_i(\beta_0,G_0)-(1-\tau)Z_if_Z(Z_i)\psi(Y_i,\delta
_i,Z_i'\beta_0,Z_i),$
and
\[
\psi(Y_i,\delta_i,t,z)= \int_0^{Y_i\wedge
t}\frac{-g_0(s|z)\,\mathrm{d}s}{\{G_0(s|z)\}^2\{1-F_0(s|z)\}}+\frac{(1-\delta
_i)I(Y_i\leq
t)}{G_0(Y_i|z)\{1-F_0(Y_i|z)\}}.
\]
\end{lemma}

\begin{pf} By the definition of $\Gamma_2$, a direct
calculation gives that
\begin{equation}\label{A4}
\Gamma_2(\beta_0,G_0)[G-G_0]=-(1-\tau) E Z\{
G(Z'\beta_0|Z)-G_0(Z'\beta_0|Z) \}/G_0(Z'\beta_0|Z).
\end{equation}

From Theorem 2.3 of Gonzalez-Manteiga and Cadarso-Suarez~\cite{G1994} and
the proof of Theorem~2 in Wang and Wang~\cite{WangWang2009}, using assumptions
C3--C7, we have that
\begin{eqnarray}\label{A5}
\hat{G}(t|z)-G_0(t|z)&=&\frac1{nh_n}\sum_{i=1}^n
K\biggl(\frac{z-Z_i}{h_n}\biggr)G_0(t|z)\psi(Y_i,\delta
_i,t,z)
\nonumber
\\[-8pt]
\\[-8pt]
\nonumber
&&{}+\mathrm{O}_p\biggl(
\biggl( \frac{\log n}{nh_n}\biggr)^{3/4}+h_n^2 \biggr).
\end{eqnarray}

Plugging (\ref{A5}) into (\ref{A4}), using standard change of variables and
Taylor expansion arguments, we obtain that $
\Gamma_2(\beta_0,G_0)[\hat{G}-G_0]=
-(1-\tau) \frac1n\sum_{i=1}^n
Z_if_Z(Z_i)\psi(Y_i,\delta_i,Z_i'\beta_0,Z_i)+\mathrm{o}_p(n^{-1/2}). $
Therefore, we have
$n^{1/2}(M_n(\beta_0,G_0)+\Gamma_2(\beta_0,G_0)[\hat{G}-G_0])=
n^{-1/2}\sum_{i=1}^nV_i+\mathrm{o}_p(1).
$ An application of the central limit theorem gives that
\[
n^{1/2}\bigl(M_n(\beta_0,G_0)+\Gamma_2(\beta_0,G_0)[\hat
{G}-G_0]\bigr)\stackrel
{d}{\to}
N(0,V).
\]
This proves the lemma.
\end{pf}

\begin{pf*}{Proof of Theorem \protect\ref{th2}} We verify the conditions in
Theorem 2 in Chen \textit{et al.}~\cite{Chen2003}. Their condition~(2.1) can be
easily verified by the subgradient condition of quantile regression
(Koenker~\cite{Koenker2005}). Their conditions (2.4), (2.5$'$) and (2.6) follow
directly from Lemma~\ref{lea1},~\ref{lea2} and~\ref{lea3}, respectively. From the
definition of $\Gamma_1$, we obtain that
\[
\Gamma_1=\Gamma_1(\beta_0,G_0)= \frac{\partial
M(\beta,G_0)}{\partial\beta}\bigg|_{\beta=\beta_0}
=-EZZ'f_0(Z'\beta_0|Z),
\]
which is negative definite by assumption C7. Thus
condition (2.2) in Chen \textit{et al.}~\cite{Chen2003} holds. By the routine
Taylor expansion,
we can
verify condition (2.3) in Chen \textit{et al.}~\cite{Chen2003}. Therefore, we
obtain that
$n^{1/2}(\hat\beta_n-\beta_0)\stackrel{d}{\to}
N(0,\Gamma_1(\beta_0,G_0)^{-1}V\Gamma_1(\beta_0,G_0)^{-1}).
$
The proof is complete.
\end{pf*}
\end{appendix}

\section*{Acknowledgements}
We thank the Editor, the AE and two
referees whose comments have helped to improve the paper substantially.
Leng's research is supported in part by NUS faculty research grants.
Tong's research was partly supported by an NSF China Zhongdian Project 11131002, NSFC (No.~10971015),
Key Project of Chinese Ministry of Education (No.~309007)
and the Fundamental Research Funds for the Central
Universities.\looseness=1

%

\printhistory

\end{document}